\documentclass[10pt]{amsart}
\usepackage{graphicx,amssymb,amsfonts,amsmath,amsthm,newlfont}
\usepackage{epsfig,url}

\usepackage{axodraw4j}
\usepackage{pstricks}
\usepackage{color}

\newtheorem{thm}{Theorem}

\theoremstyle{definition}

\theoremstyle{remark}

\newcommand{\Z}{\mathbb Z}

\newcommand{\R}{\mathbb R}

\begin{document}
\author{Oliver Schnetz}
\begin{title}[Triangle and box ladder graphs]{Evaluation of the period of a family of triangle and box ladder graphs}\end{title}
\begin{abstract}We prove that the period of a family of $n$ loop graphs with triangle and box ladders evaluates to
$\frac{4}{n}\binom{2n-2}{n-1}\zeta(2n-3)$.
\end{abstract}
\maketitle

\section{Introduction}
The period of a primitive logarithmically divergent Feynman graph $G$ is the scheme independent residue of the regularized amplitude.
It can be written in parametric space as follows. Number the edges of $G$ from $1$ to $N=2h_G$ (where $h_G$ is the number of independent loops in $G$),
and to each edge $e$ associate a variable $\alpha_e$. The period of $G$ is given by the convergent (projective) integral \cite{BEK}:
\begin{equation} \label{IZdef}
P_G=\int_{\alpha_i>0}\frac{\mathrm{d}\alpha_1\ldots\mathrm{d}\alpha_{N-1}}{\Psi_G(\alpha_1,\ldots,\alpha_{N-1},1)^2}  \in \R
\end{equation}
where $\Psi_G \in \Z[ \alpha_1, \ldots, \alpha_N]$ is the graph, or Kirchhoff, polynomial of $G$. It is defined by the formula
$$
\Psi_G = \sum_{ T \subset G} \prod_{e \notin T} \alpha_e,
$$
where the sum is over all spanning trees $T$ of $G$.

\begin{center}
\fcolorbox{white}{white}{
  \begin{picture}(336,75) (24,-38)
    \SetColor{Black}
    \SetWidth{0.8}
    \Vertex(53,33){2.8}
    \Vertex(64,-13){2.8}
    \Vertex(40,-13){2.8}
    \Vertex(76,5){2.8}
    \Vertex(32,24){2.8}
    \Vertex(71,24){2.8}
    \Vertex(28,3){2.8}
    \Arc(52,9)(24.413,145,505)
    \Arc(86.871,38.557)(34.324,-170.683,-108.465)
    \Arc(85.63,2.475)(26.026,124.203,213.79)
    \Arc(18.357,-68.443)(71.964,51.658,83.102)
    \Arc(15.214,0.679)(28.734,-26.183,54.255)
    \Arc(10.5,43)(43.661,-66.371,-13.241)
    \Arc(51,61.5)(42.5,-118.072,-61.928)
    \Arc(73.318,-34.136)(39.172,87.539,145.591)
    \Text(44,-35)[lb]{\normalsize{\Black{$\overline{Z_5}$}}}

    \Vertex(133,33){2.8}
    \Vertex(132,-16){2.8}
    \Vertex(158,9){2.8}
    \Vertex(108,9){2.8}
    \Vertex(150,26){2.8}
    \Vertex(151,-8){2.8}
    \Vertex(115,-8){2.8}
    \Vertex(116,26){2.8}
    \Arc(133,9)(24.413,145,505)
    \Arc(162.86,39.265)(30.51,-168.151,-97.268)
    \Arc(173.5,8.5)(28.504,142.125,217.875)
    \Arc(168.6,-28.567)(39.034,105.757,159.662)
    \Arc(131.231,-54.462)(50.493,66.95,109.953)
    \Arc(97.667,10)(24.333,-41.112,41.112)
    \Arc[clock](106.1,-16.8)(25.962,83.587,3.976)
    \Arc(97.4,46.567)(39.034,-74.243,-20.338)
    \Arc(133.405,57.177)(36.204,-120.555,-62.718)
    \Text(129,-35)[lb]{\normalsize{\Black{$\overline{Z_6}$}}}

    \Vertex(191,-17){2.8}
    \Vertex(207,6){2.8}
    \Vertex(214,-17){2.8}
    \Vertex(231,6){2.8}
    \Vertex(236,-17){2.8}
    \Vertex(252,6){2.8}
    \Line(191,-17)(236,-17)
    \Line(191,-17)(207,6)
    \Line(207,6)(214,-17)
    \Line(214,-17)(231,6)
    \Line(231,6)(236,-17)
    \Line(236,-17)(252,6)
    \Line(252,6)(207,6)
    \Arc[clock](221.5,-5.5)(33.534,-159.944,-339.944)
    \Text(211,-35)[lb]{\normalsize{\Black{$Z_5$}}}

    \Vertex(284,-17){2.8}
    \Vertex(299,6){2.8}
    \Vertex(308,-17){2.8}
    \Vertex(321,6){2.8}
    \Vertex(332,-17){2.8}
    \Vertex(345,6){2.8}
    \Vertex(355,-17){2.8}
    \Line(284,-17)(355,-17)
    \Line(284,-17)(299,6)
    \Line(299,6)(345,6)
    \Line(299,6)(308,-17)
    \Line(308,-17)(321,6)
    \Line(321,6)(332,-17)
    \Line(332,-17)(345,6)
    \Line(345,6)(355,-17)
    \Arc[clock](319.375,-7.6)(36.602,-165.119,-373.267)
    \Text(315,-35)[lb]{\normalsize{\Black{$Z_6$}}}
  \end{picture}
}
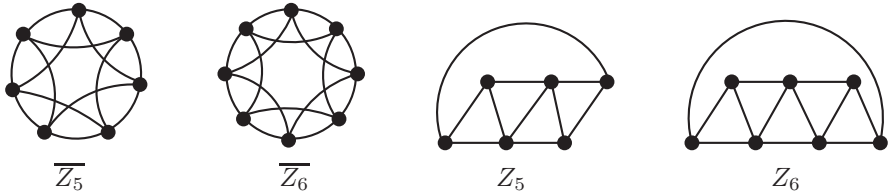
Figure 1: Completed ($\overline{Z_\bullet}$) and uncompleted ($Z_\bullet$) zig-zag graphs with 5 and 6 loops.
\end{center}
\vskip2ex

For the zig-zag graphs $Z_n$ depicted in figure 1 the periods were conjectured by D. Broadhurst and D. Kreimer \cite{BK} in 1995 as
\begin{equation} \label{IZ}
P_{Z_n} =  4 \frac{(2n-2)!}{n!(n-1)!} \left( 1 - \frac{1- (-1)^n}{2^{2n-3}}\right) \zeta(2n-3),
\end{equation}
where $\zeta(z)=\sum_{k\geq1}k^{-z}$ is the Riemann zeta function. The zig-zag conjecture was recently proved by Francis Brown and the author in \cite{ZZ}.

\begin{center}
\fcolorbox{white}{white}{
  \begin{picture}(361,286) (6,-20)
    \SetWidth{0.8}
    \SetColor{Black}
    \Vertex(10,206){2.8}
    \Vertex(30,216){2.8}
    \Vertex(30,196){2.8}
    \Vertex(10,186){2.8}
    \Vertex(10,151){2.8}
    \Vertex(30,141){2.8}
    \Vertex(10,131){2.8}
    \Vertex(30,121){2.8}
    \Vertex(50,121){2.8}
    \Vertex(70,121){2.8}
    \Vertex(50,216){2.8}
    \Vertex(70,216){2.8}
    \Vertex(120,121){2.8}
    \Vertex(140,121){2.8}
    \Vertex(160,131){2.8}
    \Vertex(140,141){2.8}
    \Vertex(160,151){2.8}
    \Vertex(120,216){2.8}
    \Vertex(140,216){2.8}
    \Vertex(95,251){2.8}
    \Vertex(160,206){2.8}
    \Vertex(140,196){2.8}
    \Vertex(160,186){2.8}
    \Line(30,216)(140,216)
    \Line(140,121)(30,121)
    \Line(10,131)(10,206)
    \Line(30,216)(30,121)
    \Line(50,216)(50,121)
    \Line(70,216)(70,121)
    \Line(120,216)(120,121)
    \Line(140,216)(140,121)
    \Line(30,216)(95,251)
    \Line(50,216)(95,251)
    \Line(70,216)(95,251)
    \Line(95,251)(120,216)
    \Line(95,251)(140,216)
    \Vertex(225,131){2.8}
    \Vertex(225,181){2.8}
    \Vertex(245,131){2.8}
    \Line(10,206)(30,216)
    \Line(10,206)(30,196)
    \Line(10,186)(30,196)
    \Line(10,151)(30,141)
    \Line(30,141)(10,131)
    \Line(10,131)(30,121)
    \Line(140,216)(160,206)
    \Line(160,206)(140,196)
    \Line(160,151)(140,141)
    \Line(140,141)(160,131)
    \Line(160,131)(140,121)
    \Line(160,206)(160,131)
    \Arc[clock](83.289,170.342)(81.504,154.055,81.739)
    \Arc(98.571,186.714)(64.385,17.43,93.18)
    \Vertex(245,181){2.8}
    \Vertex(285,131){2.8}
    \Vertex(285,181){2.8}
    \Vertex(305,181){2.8}
    \Vertex(305,131){2.8}
    \Vertex(210,146){2.8}
    \Vertex(210,131){2.8}
    \Vertex(195,161){2.8}
    \Vertex(195,146){2.8}
    \Vertex(225,211){2.8}
    \Vertex(205,206){2.8}
    \Vertex(210,221){2.8}
    \Vertex(305,211){2.8}
    \Vertex(325,206){2.8}
    \Vertex(320,221){2.8}
    \Vertex(320,131){2.8}
    \Vertex(320,146){2.8}
    \Vertex(335,146){2.8}
    \Vertex(335,161){2.8}
    \Vertex(265,251){2.8}
    \Line(225,181)(225,131)
    \Line(225,181)(210,146)
    \Line(225,181)(195,161)
    \Line(195,146)(210,146)
    \Line(210,146)(210,131)
    \Line(210,131)(320,131)
    \Line(320,131)(320,146)
    \Line(320,146)(335,146)
    \Line(335,146)(335,161)
    \Line(335,161)(305,181)
    \Line(305,181)(320,146)
    \Line(305,181)(305,131)
    \Line(225,181)(305,181)
    \Line(225,181)(225,211)
    \Line(225,211)(210,221)
    \Line(210,221)(205,206)
    \Line(205,206)(225,181)
    \Line(305,181)(305,211)
    \Line(305,211)(320,221)
    \Line(305,181)(265,251)
    \Line(285,181)(265,251)
    \Line(245,181)(265,251)
    \Line(225,181)(265,251)
    \Line(225,211)(265,251)
    \Line(210,221)(265,251)
    \Line(305,211)(265,251)
    \Line(320,221)(265,251)
    \Bezier(210,131)(145,156)(150,231)(265,251)
    \Bezier(320,131)(385,161)(385,231)(265,251)
    \Bezier(195,146)(155,181)(170,231)(265,251)
    \Bezier(335,146)(375,171)(350,231)(265,251)
    \Vertex(30,46){2.8}
    \Vertex(50,46){2.8}
    \Vertex(30,16){2.8}
    \Vertex(50,16){2.8}
    \Vertex(90,16){2.8}
    \Vertex(110,16){2.8}
    \Vertex(110,46){2.8}
    \Vertex(90,46){2.8}
    \Vertex(70,76){2.8}
    \Line(320,221)(325,206)
    \Line(325,206)(305,181)
    \Text(90,166)[lb]{\normalsize{\Black{$\ldots$}}}
    \Text(260,146)[lb]{\normalsize{\Black{$\ldots$}}}
    \Text(64,26)[lb]{\normalsize{\Black{$\ldots$}}}
    \Text(20,166)[lb]{\normalsize{\Black{$\vdots$}}}
    \Text(150,166)[lb]{\normalsize{\Black{$\vdots$}}}
    \Text(210,181)[lb]{\normalsize{\Black{$\vdots$}}}
    \Text(320,181)[lb]{\normalsize{\Black{$\vdots$}}}
    \Text(5,106)[lb]{\normalsize{\Black{$2k-1$}}}
    \Text(140,106)[lb]{\normalsize{\Black{$2\ell-1$}}}
    \Text(90,172)[lb]{\normalsize{\Black{$m$}}}
    \Text(253,151)[lb]{\normalsize{\Black{$m-1$}}}
    \Text(200,182)[lb]{\normalsize{\Black{$k$}}}
    \Text(325,182)[lb]{\normalsize{\Black{$\ell$}}}
    \Text(66,30)[lb]{\normalsize{\Black{$m$}}}
    \Text(70,100)[lb]{\normalsize{\Black{$G_{k,\ell,m}$}}}
    \Text(255,100)[lb]{\normalsize{\Black{$\widetilde{G}_{k,\ell,m}$}}}
    \Text(50,-5)[lb]{\normalsize{\Black{$\widetilde{G}_{1,1,m-1}$}}}
    \Line(140,196)(160,186)
    \Line(195,161)(195,146)
    \Line(245,181)(245,131)
    \Line(285,181)(285,131)
    \Line(30,16)(110,16)
    \Line(110,16)(110,46)
    \Line(110,46)(70,76)
    \Line(70,76)(30,46)
    \Line(70,76)(50,46)
    \Line(50,46)(50,16)
    \Line(70,76)(90,46)
    \Line(90,46)(90,16)
    \Line(30,46)(30,16)
    \Line(30,46)(110,46)
    \Arc[clock](49.457,46.868)(35.647,-123.082,-305.19)
    \Arc(90.3,46.2)(36.057,-56.883,124.263)
  \end{picture}
}
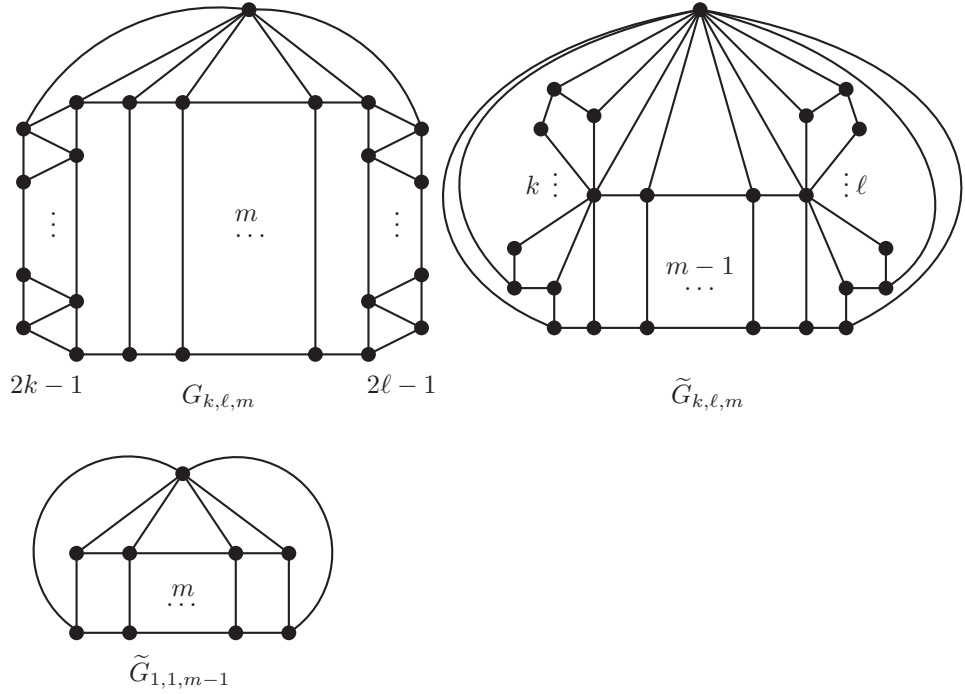
Figure 2: The $G_{k,\ell,m}$ and their planar duals $\widetilde{G}_{k,\ell,m}$ have period $\frac{4}{n}\binom{2n-2}{n-1}\zeta(2n-3)$ where $n=2(k+\ell+m)$.
\end{center}
\vskip2ex

Motivated by a conjecture in $N=4$ Super Yang-Mills theory \cite{SYM} on the period of the graph $\widehat{G}_{1,1,m-1}$ we prove that the
family $G_{k,\ell,m}$ and their planer duals $\widehat{G}_{k,\ell,m}$ have the same period as $Z_{2k+2\ell+2m}$. The graph $G_{k,\ell,m}$ depicted in figure 2 has
two vertical ladders of $2k-1$ and $2\ell-1$ triangles which are joined at their longer sides by a ladder of $m$ boxes. The three-valent vertices in the upper half are connected
to a common vertex. Its dual $\widehat{G}_{k,\ell,m}$ has two roses of $k$ and $\ell$ boxes which are joined by a ladder of $m-1$ boxes. The two-valent vertices
of the roses and the upper vertices of the box ladder are joined to a common vertex. The case $k=\ell=1$ is a ladder of $m+1$ boxes whose two-valent vertices
together with their upper three-valent vertices are connected to a common vertex.
\begin{thm}\label{main}
Let $G_{k,\ell,m}$ for $k,\ell,m\geq1$ be the family of graphs depicted in figure 2. Let $\widehat{G}_{k,\ell,m}$ be their planar duals. Then
\begin{equation}\label{maineq}
P_{G_{k,\ell,m}}=P_{\widehat{G}_{k,\ell,m}}=\frac{4}{n}\binom{2n-2}{n-1}\zeta(2n-3),
\end{equation}
where $n=2(k+\ell+m)$.
\end{thm}
\noindent The proof uses the twist-identity \cite{SchnetzCensus}.

\noindent \emph{Acknowledgements.} The author is visiting scientists at Humboldt University, Berlin.
\section{Proof of theorem \ref{main}}
\begin{center}
\fcolorbox{white}{white}{
  \begin{picture}(361,204) (6,0)
    \SetWidth{0.8}
    \SetColor{Black}
    \Vertex(10,144){2.8}
    \Vertex(30,154){2.8}
    \Vertex(30,134){2.8}
    \Vertex(10,124){2.8}
    \Vertex(10,89){2.8}
    \Vertex(30,79){2.8}
    \Vertex(10,69){2.8}
    \Vertex(30,59){2.8}
    \Vertex(50,59){2.8}
    \Vertex(70,59){2.8}
    \Vertex(50,154){2.8}
    \Vertex(70,154){2.8}
    \Vertex(120,59){2.8}
    \Vertex(140,59){2.8}
    \Vertex(160,69){2.8}
    \Vertex(140,79){2.8}
    \Vertex(160,89){2.8}
    \Vertex(120,154){2.8}
    \Vertex(140,154){2.8}
    \Vertex(95,189){2.8}
    \Vertex(160,144){2.8}
    \Vertex(140,134){2.8}
    \Vertex(160,124){2.8}
    \Line(30,154)(140,154)
    \Line(140,59)(30,59)
    \Line(10,69)(10,144)
    \Line(30,154)(30,59)
    \Line(50,154)(50,59)
    \Line(70,154)(70,59)
    \Line(120,154)(120,59)
    \Line(140,154)(140,59)
    \Line(30,154)(95,189)
    \Line(50,154)(95,189)
    \Line(70,154)(95,189)
    \Line(95,189)(120,154)
    \Line(95,189)(140,154)
    \Line(10,144)(30,154)
    \Line(10,144)(30,134)
    \Line(10,124)(30,134)
    \Line(10,89)(30,79)
    \Line(30,79)(10,69)
    \Line(10,69)(30,59)
    \Line(140,154)(160,144)
    \Line(160,144)(140,134)
    \Line(160,89)(140,79)
    \Line(140,79)(160,69)
    \Line(160,69)(140,59)
    \Line(160,144)(160,69)
    \Arc[clock](83.289,108.342)(81.504,154.055,81.739)
    \Arc(98.571,124.714)(64.385,17.43,93.18)
    \Text(75,110)[lb]{\normalsize{\Black{$\ldots$}}}
    \Text(20,104)[lb]{\normalsize{\Black{$\vdots$}}}
    \Text(150,104)[lb]{\normalsize{\Black{$\vdots$}}}
    \Line(140,134)(160,124)
    \Vertex(95,24){2.8}
    \Line(30,59)(95,24)
    \Line(50,59)(95,24)
    \Line(70,59)(95,24)
    \Line(120,59)(95,24)
    \Line(140,59)(95,24)
    \Arc(75,89)(68.007,-162.897,-72.897)
    \Arc[clock](106.324,77.088)(54.282,-8.569,-102.041)
    \Photon[double,sep=5](95,189)(95,24){2.5}{11}
    \Text(98,130)[lb]{\normalsize{\Black{$m\!-\!1$}}}
    \Text(93,199)[lb]{\normalsize{\Black{$a$}}}
    \Text(55,138)[lb]{\normalsize{\Black{$b$}}}
    \Text(55,71)[lb]{\normalsize{\Black{$0$}}}
    \Text(90,10)[lb]{\normalsize{\Black{$\infty$}}}
    \EBox(88,181)(102,195)
    \EBox(43,147)(57,161)
    \EBox(43,51)(57,67)
    \EBox(88,18)(102,31)

    \Vertex(210,144){2.8}
    \Vertex(230,154){2.8}
    \Vertex(230,134){2.8}
    \Vertex(210,124){2.8}
    \Vertex(210,89){2.8}
    \Vertex(230,79){2.8}
    \Vertex(210,69){2.8}
    \Vertex(230,59){2.8}
    \Vertex(250,59){2.8}
    \Vertex(270,59){2.8}
    \Vertex(250,154){2.8}
    \Vertex(270,154){2.8}
    \Vertex(320,59){2.8}
    \Vertex(340,59){2.8}
    \Vertex(360,69){2.8}
    \Vertex(340,79){2.8}
    \Vertex(360,89){2.8}
    \Vertex(320,154){2.8}
    \Vertex(340,154){2.8}
    \Vertex(295,189){2.8}
    \Vertex(360,144){2.8}
    \Vertex(340,134){2.8}
    \Vertex(360,124){2.8}
    \Line(230,154)(340,154)
    \Line(340,59)(230,59)
    \Line(210,69)(210,144)
    \Line(230,154)(230,59)
    \Line(270,154)(270,59)
    \Line(320,154)(320,59)
    \Line(340,154)(340,59)
    \Line(230,154)(295,189)
    \Line(250,154)(295,189)
    \Line(270,154)(295,189)
    \Line(295,189)(320,154)
    \Line(295,189)(340,154)
    \Line(210,144)(230,154)
    \Line(210,144)(230,134)
    \Line(210,124)(230,134)
    \Line(210,89)(230,79)
    \Line(230,79)(210,69)
    \Line(210,69)(230,59)
    \Line(340,154)(360,144)
    \Line(360,144)(340,134)
    \Line(360,89)(340,79)
    \Line(340,79)(360,69)
    \Line(360,69)(340,59)
    \Line(360,144)(360,69)
    \Arc(298.571,124.714)(64.385,17.43,93.18)
    \Text(275,110)[lb]{\normalsize{\Black{$\ldots$}}}
    \Text(220,104)[lb]{\normalsize{\Black{$\vdots$}}}
    \Text(350,104)[lb]{\normalsize{\Black{$\vdots$}}}
    \Line(340,134)(360,124)
    \Vertex(295,24){2.8}
    \Line(230,59)(295,24)
    \Line(250,59)(295,24)
    \Line(270,59)(295,24)
    \Line(320,59)(295,24)
    \Line(340,59)(295,24)
    \Arc[clock](306.324,77.088)(54.282,-8.569,-102.041)
    \Photon[double,sep=5](295,189)(295,24){2.5}{11}
    \Text(298,130)[lb]{\normalsize{\Black{$m\!-\!2$}}}
    \Text(293,199)[lb]{\normalsize{\Black{$a$}}}
    \Text(255,138)[lb]{\normalsize{\Black{$b$}}}
    \Text(255,71)[lb]{\normalsize{\Black{$0$}}}
    \Text(290,10)[lb]{\normalsize{\Black{$\infty$}}}
    \EBox(288,181)(302,195)
    \EBox(243,147)(257,161)
    \EBox(243,51)(257,67)
    \EBox(288,18)(302,31)
    \Arc[clock](230.618,148.99)(20.019,-168.503,-345.508)
    \Arc(230,64)(20.616,165.964,345.964)
  \end{picture}
}
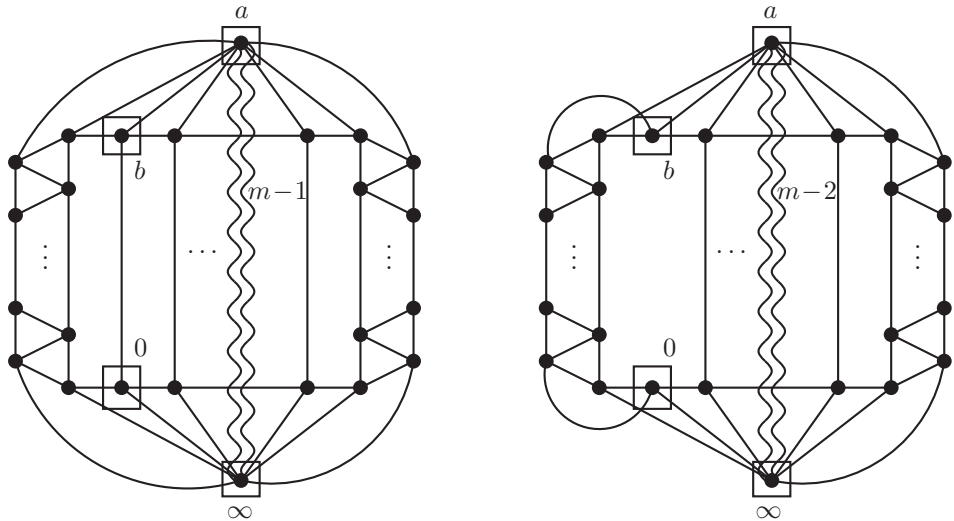
Figure 3: The completed graph $\overline{G}_{k,\ell,m}$ maps under the twist identity with respect to the boxed vertices $a,b,0,\infty$ to $\overline{G}_{k+1,\ell,m-1}$.
The curly lines symbolize propagators of negative weights $m-1$ and $m-2$, respectively.
\end{center}
\vskip2ex

\begin{proof}
First, we notice that $P_{G_{k,\ell,m}}=P_{\widehat{G}_{k,\ell,m}}$ because the period is invariant under taking planer duals \cite{BK}, \cite{SchnetzCensus}.
To prove theorem \ref{main} for $P_{G_{k,\ell,m}}$ we complete the graph by adding a vertex $\infty$ and connecting $\infty$ to all three-valent vertices.
To make the graph four-regular (make all vertices four-valent) we need to add an inverse propagator of weight $m-1$ that connects $\infty$ with $a$ in figure 3.
In position space Feynman rules a negative propagator of weight $w$ from $x$ to $y$ gives a factor $||x-y||^{2w}$ in the numerator of the integrand.
Here we need negative propagators only in intermediate steps.

We apply the twist identity \cite{SchnetzCensus} to the four vertices $a$, $b$, $0$, $\infty$ and obtain the graph on the right hand side of figure 3.
The twist identity is applied to a four vertex cut by swapping the connections of the left hand side to $a$ and $b$ and simultaneously
swapping the connections to $0$ and $\infty$. Afterwards we have to move edges to opposite sides of the four-cycle $a0b\infty$ to render the graph four-regular.
In figure 3 we had to move the edge connecting $b$ and 0 to an edge connecting $a$ and $\infty$. This new edge cancels one of the $m-1$ negative
propagators leaving a negative weight of $m-2$. After the twist we flip the left triangle ladder inside the box with vertices $b$ and $0$ and we obtain the graph
$\overline{G}_{k+1,\ell,m-1}$. Upon un-completing by removing $\infty$ we obtain
$$
P_{G_{k,\ell,m}}=P_{G_{k+1,\ell,m-1}}.
$$
By moving every second vertex in figure 1 inside the circle we see that
$$
G_{k,\ell,1}=Z_{2k+2\ell+2}.
$$
The theorem follows from (\ref{IZ}).
\end{proof}
We close this note with the remark that periods that are rational multiples of a singe zeta value are rare. The only known periods of this type are the periods
of the wheels and the zig-zags. However, with increasing loop order an increasing number of graphs can be transformed to the wheel
or the zig-zag by a sequence of twist identities and taking planar duals (the Fourier identity).


\end{document}